\documentclass[1p]{elsarticle}
\usepackage{amssymb}
\usepackage{amsfonts}
\usepackage{amsmath}

\newtheorem{theorem}{Theorem}

\newtheorem{corollary}[theorem]{Corollary}

\newtheorem{claim}[theorem]{Claim}

\newproof{pf}{Proof}

\begin{document}
\title{Asymptotically optimal bound on the adjacent vertex distinguishing edge choice number}

\author{Jakub Kwa\'sny\fnref{MNiSW}}
\ead{jkwasny@agh.edu.pl}

\author{Jakub Przyby{\l}o\fnref{grantJP,MNiSW}}
\ead{jakubprz@agh.edu.pl, phone: 048-12-617-46-38,  fax: 048-12-617-31-65}

\fntext[grantJP]{Supported by the National Science Centre, Poland, grant no. 2014/13/B/ST1/01855.}
\fntext[MNiSW]{Partly supported by the Polish Ministry of Science and Higher Education.}

\address{AGH University of Science and Technology, al. A. Mickiewicza 30, 30-059 Krakow, Poland}

\begin{abstract}
An adjacent vertex distinguishing edge colouring of a graph $G$ without isolated edges is its
proper edge colouring such that no pair of adjacent vertices meets the same set of colours in $G$.
We show that such colouring can be chosen from any set of lists associated to the edges of $G$
as long as the size of every list is at least $\Delta+C\Delta^{\frac{1}{2}}(\log\Delta)^4$, where $\Delta$ is the
maximum degree of $G$ and $C$ is a constant.
The proof is probabilistic. The same is true in the environment of total colourings.
\end{abstract}


\maketitle

\section{Introduction}

Let $G=(V,E)$ be a (simple) graph. Consider an edge colouring $c:E\to C$ where $C$ is a set of colours. For a given vertex $v\in V$, by $E(v)$ we denote the set of all edges incident with $v$ in $G$, and we set
\begin{equation}\label{palette}
S_c(v)=\{c(e):e\in E(v)\}.
\end{equation}
We shall also write $S(v)$ instead of $S_c(v)$ provided this causes no ambiguities,
and we shall call such set a \emph{palette} of $v$ or simply refer to it as a \emph{set of colours incident with} $v$.
The colouring $c$ is called \emph{adjacent vertex distinguishing} if it is proper and $S(u)\neq S(v)$ for every edge $uv\in E$.
It exists if only $G$ contains no isolated edges.
The least number of colours in $C$ necessary to provide such a colouring is
denoted by $\chi'_{a}(G)$ and called the \emph{adjacent vertex
distinguishing edge chromatic number} of $G$, see~\cite{Hatami},
and~\cite{Akbari,BalGLS,FlandrinMPSW,Zhang} for alternative notations used.
It was conjectured~\cite{Zhang} that $\chi'_{a}(G)\leq \Delta(G)+2$ for every connected graph $G$ of order at
least three different from the cycle $C_5$.
This was, e.g., positively verified by Balister et al.~\cite{BalGLS} for bipartite graphs and for graphs of maximum degree $3$, while Greenhill and Ruci\'nski proved it asymptotically almost surely for random $4$-regular graphs, see~\cite{Rucinski_regular},
and~\cite{BonamyEtAl,BuLihWang,HocqMont0,HocqMont,Hornak_planar,LiZhangChenSun,WangWang} for results concerning other particular graph classes.
In general it is known that $\chi'_a(G)\leq 3\Delta(G)$,~\cite{Akbari},
and $\chi'_a(G)\leq \Delta(G)+O(\log \chi(G))$,~\cite{BalGLS}, for every graph $G$ with no isolated edges.
Finally, Hatami~\cite{Hatami} proved the postulated upper bound up to an additive constant by showing that $\chi'_a(G)\leq \Delta+300$
for every graph $G$ with no isolated edges and with maximum degree $\Delta>10^{20}$.

Suppose now that every edge $e\in E$ is endowed with a list of available colours $L_e$.
Analogously as in the case of the classical choosability of graphs, introduced for vertex colourings by Vizing~\cite{Vizing_list}
and independently by Erd\H{o}s, Rubin and Taylor~\cite{ErdosRubinTaylor},
we define the \emph{adjacent vertex distinguishing edge choice number} of a graph $G$ (without isolated edges)
as the least $k$ so that for every set of lists of size $k$ associated to the edges of $G$
we are able to choose colours from the respective lists to obtain an adjacent vertex distinguishing edge colouring of $G$.
We denote it by ${\rm ch}'_a(G)$.
This was already investigated under different notations e.g. in~\cite{HornakWozniak-list_avd} (for a few simple graph classes) and in~\cite{Przybylo_CN_1,Przybylo_CN_2}, where it was proved that ${\rm ch}'_a(G)\leq 2\Delta(G)+{\rm col}(G)-1$ and ${\rm ch}'_a(G)\leq \Delta(G)+3{\rm col}(G)-4$ (where ${\rm col}(G)$ denotes the colouring number of $G$, hence
${\rm col}(G)\leq\Delta(G)+1$) for every graph $G$ without isolated edges.
The latter of these results implies an upper bound of the form ${\rm ch}'_a(G)\leq\Delta(G)+K$
with a constant $K$ for many classes of graph, e.g. for planar graphs.
In fact in~\cite{HornakWozniak-list_avd} it was boldly conjectured that we always have ${\rm ch}'_a(G)=\chi'_a(G)$.
This refers to the well known List Colouring Conjecture, posed independently by several  researchers  (see~\cite{JensenToft_Colorings},
Section 12.20), that we always have ${\rm ch}'(G)=\chi'(G)$ where ${\rm ch}'(G)$ denotes the edge choosability (or equivalently the list chromatic index or the edge choice number) of $G$.
If proved, the List Colouring Conjecture combined with Vizing's Theorem would imply that ${\rm ch}'(G)\leq \Delta(G)+1$.
So far the best result concerning such supposed upper bound is due to Molloy and Reed~\cite{molloy-reed-nolc}.
It implies that for every graph $G$ with maximum degree $\Delta$,
${\rm ch}'(G)=\Delta+O(\Delta^{\frac{1}{2}}(\log\Delta)^4)$ where by `$\log$' we mean the natural logarithm in this paper.
Note that by definition, ${\rm ch}'_a(G)\geq{\rm ch}'(G)$.
We shall prove that an upper bound of the same form as the one in~\cite{molloy-reed-nolc} above is valid also in the case of ${\rm ch}'_a(G)$, see Theorem~\ref{main_theorem_list_Hatami} below.

Distinguishing by colour palettes was also considered for general, not necessarily proper edge colourings.
In such a setting we strive at distinguish adjacent vertices by their respective \emph{multisets} of incident colours.
It is believed that three colours are always sufficient for this goal for any graph without isolated edges, while the best result thus far is four, see~\cite{Aldred}.
This problem is deeply related and was first considered in context of the so-called \emph{1--2--3 Conjecture}, where we use integer colours and search for a minimum $k$ so that a (not necessarily proper) edge colouring $c:E\to\{1,2,\ldots,k\}$ exists such that adjacent vertices in a graph meet distinct sums of their incident colours, see~\cite{123KLT} for details. It was conjectured there that
integers $1$, $2$ and $3$ are sufficient for any graph without isolated edges for this aim.
Thus far it is known that the labels $1,2,3,4,5$ always suffice, see~\cite{KalKarPf_123}
(while the mentioned result for multisets implies that the same can be achieved using a set of four
\emph{real} labels -- it is sufficient to choose these to be independent over the field of rational numbers).
Surprisingly, no constant upper bound is known in a natural list correspondent of this concept, though it is believed that lists (of real numbers) of cardinality three should always suffice, see~\cite{BarGrNiw}.
This is on the other hand known to hold for a total analogue
of such a choosability problem, see~\cite{WongZhu23Choos} for details.

The common source of all these problems is the graph invariant called the \emph{irregularity strength} of a graph $G$, i.e. the least $k$ so that a colouring $c:E\to\{1,2,\ldots,k\}$ exists attributing every vertex of $G$ a distinct sum of its incident colours, or equivalently the least $k$ such that we are able to multiply some of the edges of a given graph $G$ -- each at most $k$ times in order to construct of $G$ an \emph{irregular multigraph}, i.e. a multigraph with pairwise distinct vertex degrees, see e.g.~\cite{Aigner,Chartrand,Lazebnik,Faudree2,Frieze,KalKarPf,MajerskiPrzybylo2,Nierhoff}.
This concept was motivated and stemmed from the fact that there are no irregular \emph{graphs} at all, except the trivial one--vertex case and the related research of Chartrand, Erd\H{o}s and Oellermann~\cite{ChartrandErdosOellermann} on possible alternative definitions of irregularity in graphs.

\section{Preliminaries}
Our main result is the following.
\begin{theorem}\label{main_theorem_list_Hatami}
There is a constant $C$ such that
$${\rm ch}'_a(G)\leq \Delta + C\Delta^\frac{1}{2}(\log\Delta)^4$$
for every graph $G$ with maximum degree $\Delta$ and without isolated edges.
\end{theorem}
To prove this we shall apply several times probabilistic approach.
As a starting point of our construction of a desired colouring of a given graph we shall however need the following result of Molloy and Reed
concerning the more general setting of hypergraphs. There the \emph{maximum codegree} of a hypergraph denotes the greatest number of its edges containing the same pair of vertices.
\begin{theorem}[\cite{molloy-reed-nolc}]\label{MolloyReedThChoosabilityCodegree}
For all $k$, there is a constant $C_k$ (depending on $k$) such that any $k$--uniform hypergraph of maximum degree $\Delta$ and maximum codegree $B$ has list chromatic index at most $\Delta+C_kB^\frac{1}{k}\Delta^{1-\frac{1}{k}}\left(\log\frac{\Delta}{B}\right)^4$.
\end{theorem}
Note that this theorem implies the mentioned above best known (in terms of $\Delta$) upper bound for
the lengths of lists associated to the edges of a given graph $G$ from which one may always properly edge colour the graph, i.e. $\Delta+O(\Delta^\frac{1}{2}(\log\Delta)^4)$.
We wish to show the same size of lists (up to the multiplicative constant at the second order term) is also sufficient to additionally distinguish adjacent vertices with their corresponding colour palettes.

We shall use Theorem~\ref{MolloyReedThChoosabilityCodegree} above with $k=2$ and $B\leq 2$, hence for the case of multigraphs (or graphs) with edge multiplicity at most two\footnote{Though in~\cite{molloy-reed-nolc} the family of edges of a hypergraph was defined as a set (not multiset), i.e. excluding existence of two copies of the same edge, it can be verified that the thesis of Theorem~\ref{MolloyReedThChoosabilityCodegree} holds by the same argument if we admit multiple hyperedges.}. This shall be useful, as our Theorem~\ref{main_theorem_list_Hatami} is in fact directly implied by the following one, which we shall prove below, first discussing the mentioned implication in Section~\ref{SectionMultigraphGraph}. (The parameter ${\rm ch}'_a(G)$ for a \emph{multigraph} $G$ is defined the same way as for a graph, by means of~(\ref{palette}), where $E(v)$ denotes the set of all edges incident with $v$ -- counting in all edges joining $v$ with the same neighbour).
\begin{theorem}\label{main_theorem_list_Hatami_auxiliary}
There is a constant $C_0$ and $\Delta_0$ such that
$${\rm ch}'_a(G)\leq \Delta + C_0\Delta^\frac{1}{2}(\log\Delta)^4$$
for every multigraph $G$ with edge multiplicity at most $2$, maximum degree $\Delta\geq\Delta_0$ and minimum degree $\delta\geq\frac{\Delta}{4}$.
\end{theorem}

In our approach, similarly as in~\cite{Hatami} we start from any proper edge colouring,
and then recolour in a few steps some part of the edges, but in a different way than in that paper
(not to mention that we must choose colours from specified and potentially diversified lists).
We shall proceed in four stages as follows:
\begin{enumerate}
 \item[(I)]
     We first reserve some portion of colours in each given edge list for a later use. At the same time we remove some additional colours from each of these so that for any partial colouring $c$ from the leftovers of the lists (still including great majority of all original colours) and any its supplementation $c'$ from the respective sets of
     reserved colours, there can be no colour conflict between adjacent edges coloured under $c$ and $c'$.
 \item[(II)] We then fix by Theorem~\ref{MolloyReedThChoosabilityCodegree} any proper colouring from the leftovers of the lists and randomly uncolour some small portion of the edges, however large enough so that afterwards adjacent vertices are not only distinguished, but the symmetric difference of their partial colour palettes  is relatively large.
 \item[(III)]
      Next, we use the naive colouring procedure on the uncoloured edges, randomly choosing colours from their respective reserved lists and uncolouring adjacent edges coloured the same. We guarantee
      that afterwards the symmetric differences between adjacent partial colour palettes are large
      compared to vertex degrees in the uncoloured subgraph.
 \item[(IV)] At the end, any greedy choice of the reserved colours for still uncoloured edges satisfies our target requirement, 
      as long as we choose them so that the obtained edge colouring is proper.
\end{enumerate}

In order to control our random process we shall use a few tools of the probabilistic method,
the Lov\'asz Local Lemma, see e.g.~\cite{AlonSpencer},
the Chernoff Bound, see e.g.~\cite{JansonLuczakRucinski} (Th. 2.1, page 26)
and Talagrand's Inequality, see e.g.~\cite{MolloyReed_GoodTalagrand}:
\begin{theorem}[\textbf{The Local Lemma}]
\label{LLL-symmetric}
Let $A_1,A_2,\ldots,A_n$ be events in an arbitrary pro\-ba\-bi\-li\-ty space.
Suppose that each event $A_i$ is mutually independent of a set of all the other
events $A_j$ but at most $D$, and that ${\mathbf Pr}(A_i)\leq p$ for all $1\leq i \leq n$. If
$$ ep(D+1) \leq 1,$$
then $ {\mathbf Pr}\left(\bigcap_{i=1}^n\overline{A_i}\right)>0$.
\end{theorem}
\begin{theorem}[\textbf{Chernoff Bound}]\label{ChernofBoundTh}
For any $0\leq t\leq np$,
$${\mathbf Pr}({\rm BIN}(n,p)>np+t)<e^{-\frac{t^2}{3np}}~~{and}~~{\mathbf Pr}({\rm BIN}(n,p)<np-t)<e^{-\frac{t^2}{2np}}\leq e^{-\frac{t^2}{3np}}$$
where ${\rm BIN}(n,p)$ is the sum of $n$ independent Bernoulli variables, each equal to $1$ with probability $p$ and $0$ otherwise.
\end{theorem}

\begin{theorem}[\textbf{Talagrand's Inequality}]\label{TalagrandsInequalityTotal}
Let $X$ be a non-negative random variable determined by $l$ independent trials $T_1,\ldots,T_l$.
Suppose there exist constants $c,k>0$ such that for every set of possible outcomes of the trials, we have:
\begin{itemize}
\item[1.] changing the outcome of any one trial can affect $X$ by at most $c$, and
\item[2.] for each $s>0$, if $X\geq s$ then there is a set of at most $ks$ trials whose outcomes certify that $X\geq s$.
\end{itemize}
Then for any $t\geq 0$, we have
\begin{equation}\label{TalagrandsInequality}
{\mathbf Pr}(|X-{\mathbf E}(X)|>t+20c\sqrt{k{\mathbf E}(X)}+64c^2k)\leq 4e^{-\frac{t^2}{8c^2k({\mathbf E}(X)+t)}}.
\end{equation}
\end{theorem}

\section{Colouring multigraphs yields graph colourings}\label{SectionMultigraphGraph}

Suppose Theorem~\ref{main_theorem_list_Hatami_auxiliary} holds, and let $C_0$ and $\Delta_0$ be the constants from it.
Set $C:=\max\{C_0,3\Delta_0,10\}$.
Let $G=(V,E)$ be a graph without isolated edges and let $\left\{L_e\right\}_{e\in E}$ be a set of lists of lengths $\Delta+\lfloor C\Delta^\frac{1}{2}(\log\Delta)^4\rfloor$, where $\Delta$ is the maximum degree of $G$ ($\Delta\geq 2$).
We show there exists an adjacent vertex distinguishing edge colouring of $G$ from these lists.
If $\Delta< \Delta_0$, we colour the graph greedily component by component.
More precisely, for every consecutive edge $uv$ of a component we choose a colour from $L_{uv}$ distinct from all the colours
already assigned to its at most $2(\Delta-1)$ adjacent edges such that the resulting partial palette at $u$ (i.e. the set of all colours already assigned to incident edges of $u$) is distinct from the (partial) palettes of all its neighbours except possibly $v$ and so that the palette of $v$ is distinct from the palettes of its neighbours except possibly $u$. This is feasible, as $|L_{uv}|> 4(\Delta-1)$. As $G$ has no components of size $1$, at the end we obtain an adjacent vertex distinguishing edge colouring of $G$.

So assume that $\Delta\geq \Delta_0$. For every vertex $u$ of $G$ with $d(u)<\Delta/4$ which has exactly one neighbour $v$ with $d(v)<\Delta/4$ we contract the edge $uv$ (keeping the possible multiple edges). The resulting multigraph $G'$ has maximum degree $\Delta$ and edge multiplicity at most $2$.
In the remaining part of the paper we shall abbreviate the term `multigraph' and usually write `graph' instead.
If the minimum degree of $G'$ is less than $\Delta/4$ we take two copies of $G'$ and add an edge between every vertex of degree less than $\Delta/4$ in the first copy and its analogue in the second copy. We repeat such procedure (taking again two copies of the resulting graph) if necessary until we obtain a graph $G''$ with minimum degree $\delta\geq\Delta/4$. We associate any list of colours of lengths $\Delta+\lfloor C\Delta^\frac{1}{2}(\log\Delta)^4\rfloor$ to edges without such lists assigned yet. By Theorem~\ref{main_theorem_list_Hatami_auxiliary} there is an adjacent vertex distinguishing edge colouring of $G''$ from the given lists. This restricted to the edges of $G'$ yields its proper edge colouring $c$ from the original lists under which every vertex of degree at least $\Delta/4$ is set distinguished (has a distinct palette) from all its neighbours.
This remains valid also in $G$ if we uncontract the previously contracted edges and uncolour all edges of $G$ with both ends in the set $S=\{v\in V:d(v)<\Delta/4\}$. Let $H=(S,E')$ be the graph induced by the set $E'$ of all the uncoloured edges in $G$.
Note that the only components of $H$ with size $1$ are formed by the previously contracted edges $uv$, and thus such $u$ and $v$ are (and shall remain) set distinguished in $G$ as their partial palettes are at this point disjoint. Therefore we may analogously as above greedily choose new colours for the edges of $H$ from their respective lists so that afterwards also all vertices in $S$
are set distinguished from their neighbours in $G$. (Note in particular that as we only need to distinguish neighbours of the same degree in $G$, any vertex from $S$ forming a $1$-vertex component of $H$ is trivially set distinguished from all its neighbours in $G$.)

\section{Proof of Theorem~\ref{main_theorem_list_Hatami_auxiliary}}

We do not specify $\Delta_0$. We simply assume throughout the proof that it is a large enough constant that all explicit inequalities below hold; in particular $\Delta_0\geq 9$. Let $C'=C_2 \sqrt{2}$ where $C_2$ is a constant from Theorem~\ref{MolloyReedThChoosabilityCodegree}
(exploited further with $k=2$ and $B=2$) and set $C_0=C'+4$.
Let $G=(V,E)$ be a (multi)graph of codegree at most $2$ with maximum degree $\Delta\geq \Delta_0$ and minimum degree $\delta\geq \Delta/4$ (hence without a component of order less than $3$)
endowed with a set $\left\{L_e\right\}_{e\in E}$ of lists, each of cardinality $\Delta+\lfloor C_0\Delta^\frac{1}{2}(\log\Delta)^4\rfloor$.
Further on it shall be convenient for us to denote a given edge e.g. by $uv$, bearing in mind that there might be two different edges joining $u$ and $v$ in $G$.

\subsection{Step I: Reserving colours}

For each $v\in V$ we choose a set $R_v\subset \bigcup_{e\in E(v)} L_{e}$
independently placing in it every colour from $\bigcup_{e\in E(v)} L_{e}$
with probability $p_1 = \frac{(\log\Delta)^4}{\Delta^\frac{1}{2}}$.
For every edge $uv\in E$ we then define the \emph{list of reserved colours} as $R_{uv} = R_u \cap R_v\cap L_{uv}$ and the \emph{leftover list} as $L'_{uv} = L_{uv} \smallsetminus (R_u \cup R_v)$.

\begin{claim}\label{cl: reserve} The sets $R_w$, $w\in V$ can be chosen so that for every $e\in E$:
\begin{itemize}
 \item[(i)] $|R_e|\geq \frac12 (\log\Delta)^8$;
 \item[(ii)] $|L'_e|\geq \Delta + C'\Delta^\frac{1}{2}(\log\Delta)^4$.
\end{itemize} \end{claim}
\begin{pf} For each edge $uv\in E$, the random variable $X_{uv}=|R_u\cap R_v \cap L_{uv}|$ is a sum of $\Delta+\lfloor C_0\Delta^\frac{1}{2}(\log\Delta)^4\rfloor$ independent 0--1 random variables (corresponding to the elements of $L_{uv}$)
with probability of picking $1$ equal to $p^2_1$.
By the Chernoff Bound we thus obtain that:
\begin{equation}\label{IneqXuv}
{\mathbf Pr}\left(X_{uv}<\frac12 (\log\Delta)^8\right) \leq {\mathbf Pr}\left(X_{uv}<\frac12 (\log\Delta)^8+1\right)
< e^{-\frac{(\log\Delta)^8}{5}}\leq\frac{1}{\Delta^2}
\end{equation}
(for $\Delta$ sufficiently large).

Now for any edge $uv\in E$ we denote a random variable $Y_{uv}=|(R_u\cup R_v) \cap L_{uv}|$. As above, it can be regarded as a sum of $\Delta+\lfloor C_0\Delta^\frac{1}{2}(\log\Delta)^4\rfloor$ independent Bernoulli trials with the probability of success equal to $(2p_1-p^2_1)$. As thus ${\mathbf E}(Y_{uv})=2\Delta^\frac12(\log\Delta)^4+\Theta((\log\Delta)^8)$, again by the Chernoff Bound we obtain:
\begin{equation}\label{IneqYuv}
{\mathbf Pr}\left(Y_{uv}>4\Delta^\frac12(\log\Delta)^4-1\right)<e^{-\frac{\Delta^\frac12(\log\Delta)^4}{2}}\leq\frac{1}{\Delta^2}.
\end{equation}

Consider all events of the form $X_e<\frac12 (\log\Delta)^8$ and $Y_{e}>4\Delta^\frac12(\log\Delta)^4-1$ with $e\in E$.
Note that any of these is mutually independent of all such events but the ones associated with edges $f$ such that $e\cap f\neq \emptyset$,
i.e. all other events but at most $4\Delta$. Since $e\frac{1}{\Delta^2}(4\Delta+1)\leq 1$, by~(\ref{IneqXuv}), (\ref{IneqYuv}) and the Lov\'asz Local Lemma, with positive probability none of these events occurs. The thesis follows by our choice of $C_0$.
$\Box$
\end{pf}

\subsection{Step II: Uncolouring edges}

Let $R_e$ and $L'_e$, $e\in E$ be sets complying with Claim~\ref{cl: reserve}.
By the choice of $C'$ and Theorem~\ref{MolloyReedThChoosabilityCodegree} there exists
a proper edge colouring of $G$ from the lists $L'_e$, $e\in E$. We fix any such colouring
and then randomly and independently uncolour each edge of $G$ with probability $p_2 = \frac{(\log\Delta)^2}{\Delta}$. For every $v\in V$, let $U_v$ be the set of uncoloured edges incident with $v$ in $G$. By $S_v$ in turn we denote the obtained \emph{partial palette} of $v$, i.e. the set of colours of the edges in $E(v)\smallsetminus U_v$. Let $S_1\triangle S_2:=(S_1\cup S_2)\smallsetminus (S_1\cap S_2)$ denote the \emph{symmetric difference} of any sets $S_1,S_2$.

\begin{claim} \label{cl: step2} We can uncolour the edges of $G$ so that:
\begin{itemize}
 \item[(iii)] $\frac18 (\log\Delta)^2 \leq |U_v| \leq \frac32 (\log\Delta)^2$ for every $v\in V$;
 \item[(iv)] $|S_u \triangle S_v| \geq \frac1{16} (\log\Delta)^2$ for every $uv\in E$ with $d(u)=d(v)$.
\end{itemize} \end{claim}
\begin{pf} For each $v$, the random variable $|U_v|$ is a sum of at least $\frac{\Delta}4$ and at most $\Delta$ independent binary random variables equal to 1 with probability $p_2 = \frac{(\log\Delta)^2}{\Delta}$, therefore $\frac14 (\log\Delta)^2 \leq {\mathbf E}(|U_v|) \leq (\log\Delta)^2$. The Chernoff Bound thus yields:
\begin{equation}\label{IneqUv}
{\mathbf Pr}\left(|U_v|<\frac18 (\log\Delta)^2\right)+{\mathbf Pr}\left(|U_v|>\frac32 (\log\Delta)^2\right)<e^{-\frac{(\log\Delta)^2}{32}} +e^{-\frac{(\log\Delta)^2}{48}} \leq \frac{1}{\Delta^3}.
\end{equation}
Let for every $v\in V$, $A_v$ denote the event that $(iii)$ does not hold for $v$. By~(\ref{IneqUv}) above, ${\mathbf Pr}(A_v)<\Delta^{-3}$.

Now we shall temporarily need an orientation of the edges.
Fix any ordering $v_1,\dots,v_l$ of $V$ and orient every edge from its end with lower index to the end with higher index.
For every edge $uv\in E$ with $d(u)=d(v)$ and oriented, without loss of generality, from $u$ to $v$ we denote by $B_{\overrightarrow{uv}}$ the event that $(iv)$ does not hold for $uv$.
Since we uncolour the edges independently, while estimating the probability of a given event $B_{\overrightarrow{uv}}$ we may assume that the uncolourings are being committed in any order convenient for our analysis.
We shall in fact bound the following probability:
\begin{equation}\label{ConditionalInequality}
{\mathbf Pr}(B_{\overrightarrow{uv}}\wedge \overline{A_v}) \leq {\mathbf Pr}\left(|S_u \triangle S_v| < \frac1{16} (\log\Delta)^2 \bigg| \frac18 (\log\Delta)^2 \leq |U_v| \leq \frac32 (\log\Delta)^2\right).
\end{equation}
The analysis of the latter of these is easier to follow if one assumes the edges incident with $v$ have been (randomly) uncoloured first, and only just then we draw the edges incident with $u$ (except one or two $uv$) to be uncoloured. (The same can be formalized by means of the law of total probability.)
By the conditional assumption, $|U_v| \geq \frac18 (\log\Delta)^2$, and thus, just before uncolouring the remaining (except at most two $uv$) edges incident with $u$, there are at least $\frac18 (\log\Delta)^2-2$ colours
associated to the edges incident with $u$ (other than $uv$) that do not appear on the edges incident with $v$.
Denote the set of these colours by $R$ and set $r=|R|$, hence $r\geq\frac18 (\log\Delta)^2-2$.
In order to have $|S_u \triangle S_v| < \frac1{16} (\log\Delta)^2$ satisfied, we thus need to uncolour at least
$\left\lceil r-\frac1{16} (\log\Delta)^2 \right\rceil$ out of $r$ edges coloured with the elements of $R$.
Therefore, by~(\ref{ConditionalInequality}),
\begin{equation}\label{ConditionalInequality2}
{\mathbf Pr}(B_{\overrightarrow{uv}}\wedge \overline{A_v}) \leq \max_{r\ge \frac18 (\log\Delta)^2-2} \binom{r}{\left\lceil r-\frac1{16} (\log\Delta)^2 \right\rceil} p_2^{\left\lceil r-\frac1{16} (\log\Delta)^2 \right\rceil}.
\end{equation}
Set $A= \left\lfloor \frac1{16} (\log\Delta)^2 \right\rfloor$.
In order to upper-bound the above it is thus sufficient to
maximize a function $f(r)=\binom{r}{r-A} p_2^{r-A}$ for integers $r\ge 2A-2$, for which we have:
$$\frac{f(r+1)}{f(r)} = \frac{(r+1)!A! (r-A)! p_2^{r+1-A}}{r!A!(r+1-A)! p_2^{r-A}} =  \frac{(r+1)p_2}{(r+1-A)} < 1$$
(for $\Delta$ sufficiently large). Therefore the function $f(r)$ is decreasing for integers $r\ge 2A-2$.
By~(\ref{ConditionalInequality2}) and the well-known fact that $\binom{2n}{n} \leq 4^n$ for $n\geq 1$, we thus obtain that:
\begin{equation}\label{ConditionalInequality3}
{\mathbf Pr}(B_{\overrightarrow{uv}}\wedge \overline{A_v}) \leq  \binom{2A-2}{A-2} p_2^{A-2} \leq  \binom{2(A-3)}{A-3} p_2^{A-3} \leq (4p_2)^{A-3} < \frac{1}{\Delta^3}.
\end{equation}
Now, every event $A_w$ and $(B_{\overrightarrow{uv}} \wedge \overline{A_v})$ is mutually independent of all other such events except
those associated with a vertex (and possibly some other vertex in case of events of the second type) at distance at most $1$ from $w$ or, resp. $u$ or $v$, i.e., all except at most $4\Delta^2$ other events. By~(\ref{IneqUv}) and~(\ref{ConditionalInequality3}) and the Local Lemma we thus conclude we may uncolour some edges of $G$ so that none of these events (hence also none of $A_w$ and $B_{\overrightarrow{uv}}$) holds -- the thesis follows.
$\Box$
\end{pf}

\subsection{Step III: Colouring with the reserved colours}

Let $U$ be a set of all uncoloured edges in $G$ consistent with Claim~\ref{cl: step2} above.
Independently for every edge $e\in U$ we now randomly and with equal probability choose a colour from its list of reserved colours $R_e$.
After finishing this process we uncolour every edge (in $U$) coloured the same as any of its adjacent edges (in $U$).
Note that consequently, by our choice of the lists of reserved colours, the obtained partial edge colouring of $G$ shall be proper.
We denote by $U'$ the set of uncoloured edges after this procedure. We also denote by $U'_v$ the set of edges in
$U'$ incident with a vertex $v$, while $S'_v$ shall denote the partial colour palette of $v$ at the end of this step.

\begin{claim}\label{c2: step3}
It is possible to choose colours from the reserved lists so that:
\begin{itemize}
 \item[(v)] $|U'_v| \leq \frac14\cdot\frac1{32} (\log\Delta)^2$ for every $v\in V$;
 \item[(vi)] $|S'_u \triangle S'_v| \geq \frac1{32} (\log\Delta)^2$ for every $uv\in E$ with $d(u)=d(v)$.
\end{itemize}
\end{claim}

\begin{pf}
For any given vertex $v\in V$, let $A'_v$ denote the event that $|U'_v| > \frac1{128} (\log\Delta)^2$.
Note that as $|R_e|\geq \frac12(\log\Delta)^8$ for every edge $e$ by Claim~\ref{cl: reserve} $(i)$, the probability that an edge $e\in U_v$ shall be uncoloured at this stage equals at most $|U_v|\cdot 2(\log\Delta)^{-8}$. Therefore, since $|U_v| \leq \frac32 (\log\Delta)^2$ by Claim~\ref{cl: step2} $(iii)$,
\begin{align*} {\mathbf E}(|U'_v|) \leq |U_v|^2 \frac2{(\log\Delta)^8} \leq \frac9{2 (\log\Delta)^{4}}\le \frac13 \cdot \frac1{128} (\log\Delta)^2. \end{align*}
In order to bound the probability of $A'_v$ we shall use Talagrand's Inequality.
Note that the random variable $|U'_v|$ is determined by at most $\Delta^2$ independent trials, i.e. the choices of the reserved colours for the edges in $U_v$ and their adjacent edges at this step. A change of the result of any such trial may alter $|U'_v|$ by at most $2$, while the fact that $|U'_v|\geq s$ can be certified by the outcomes of at most $2s$ trials (each corresponding to one uncoloured edge incident with $v$ or an adjacent edge of such an edge with the same colour assigned).
By Theorem~\ref{TalagrandsInequalityTotal} we thus obtain:

\begin{eqnarray} {\mathbf Pr} (A'_v) &\leq&
{\mathbf Pr} \left( |U'_v|> \frac13 \cdot \frac1{128} (\log\Delta)^2 + \frac13 \cdot \frac1{128} (\log\Delta)^2 + 40\sqrt{2\cdot \frac13 \cdot \frac1{128} (\log\Delta)^2} + 512 \right) \nonumber\\
&\leq& {\mathbf Pr} \left( |U'_v|> {\mathbf E}(|U'_v|) + \frac13 \cdot \frac1{128} (\log\Delta)^2 + 20\cdot 2\sqrt{2 {\mathbf E}(|U'_v|)} + 64\cdot 2^2\cdot 2 \right)\nonumber\\
&\leq& 4e^{-\frac{\left(\frac13 \cdot \frac1{128} (\log\Delta)^2\right)^2}{64\left({\mathbf E}(|U'_v|)+\frac13 \cdot \frac1{128} (\log\Delta)^2\right)}}
\leq
4e^{-\frac{\left(\frac13 \cdot \frac1{128} (\log\Delta)^2\right)^2}{64\left(2\cdot\frac13 \cdot \frac1{128} (\log\Delta)^2\right)}} < \frac{1}{\Delta^5}.
\label{PrU'bound}
\end{eqnarray}

For any edge $uv\in E$ with $d(u)=d(v)$, let $B'_{uv}$ denote the event that $|S'_u \triangle S'_v| < \frac1{32} (\log\Delta)^2$.
By Claim~\ref{cl: step2}~$(iv)$, $|S_u \triangle S_v| \geq \frac1{16} (\log\Delta)^2$. We fix any subset $S_{uv}\subseteq (S_u \triangle S_v)$ with
$|S_{uv}| = \lfloor\frac1{16} (\log\Delta)^2\rfloor$.
Note that $B'_{uv}$ implies that at least $\lfloor\frac1{32} (\log\Delta)^2\rfloor$ colours from $S_{uv}$ have been assigned to edges in $U_u\cup U_v$ at this step
(despite the fact that some of these could also be uncoloured later). The necessary condition for this
is that all edges of some subset of $\lfloor\frac1{32} (\log\Delta)^2\rfloor$ elements of $U_u\cup U_v$ have been assigned colours from $S_{uv}$.
As by Claim~\ref{cl: step2} $(iii)$, $|U_u\cup U_v| \leq 3 (\log\Delta)^2$ while by Claim~\ref{cl: reserve} $(i)$, $|R_e|\geq \frac12(\log\Delta)^8$ for every edge $e$, we obtain that:

\begin{eqnarray} {\mathbf Pr} \left(B'_{uv}\right) &\leq& \binom{\left\lfloor 3 (\log\Delta)^2\right\rfloor}{\left\lfloor\frac1{32} (\log\Delta)^2\right\rfloor} \left(\frac{\frac1{16}(\log\Delta)^2}{\frac12 (\log\Delta)^{8}}\right)^{\left\lfloor\frac1{32} (\log\Delta)^2\right\rfloor}\nonumber\\
&\leq& \left(3 (\log\Delta)^2 \right)^{\left\lfloor\frac1{32} (\log\Delta)^2\right\rfloor} \cdot \left(\frac1{8 (\log\Delta)^{6}}\right)^{\left\lfloor\frac1{32} (\log\Delta)^2\right\rfloor} \nonumber\\
& =& \left(\frac{3}{8(\log\Delta)^{4}}\right)^{\left\lfloor\frac1{32} (\log\Delta)^2\right\rfloor} < \frac1{\Delta^5}.\label{PrB'uv}
\end{eqnarray}

Analogously as previously, every event $A'_w$ and $B'_{uv}$ is mutually independent of all other events except
those associated with a vertex (and possibly some other vertex in case of events of the second type) at distance at most $3$ from $w$ or, resp. $u$ or $v$, i.e., all except at most $4\Delta^4$ other events. By~(\ref{PrU'bound}) and~(\ref{PrB'uv}) and the Local Lemma we thus obtain the thesis.
$\Box$
\end{pf}

\subsection{Step IV: Colouring the remaining edges}

Fix any (partial) edge colouring of $G$ consistent with the thesis of Claim~\ref{c2: step3}.
Now we may colour the remaining edges greedily. To do this we analyze one after another each edge $e=uv$ in $U'$ and colour it with any colour from $R_e$
which is not used by some of its adjacent edges. This is feasible, as by our construction and the choice of the lists of reserved colours, a colour from $R_e$ cannot be assigned to any edge in
$(E(u)\cup E(v))\smallsetminus (U_u\cup U_v)$, while $|U_u\cup U_v| \leq 3 (\log\Delta)^2$ by Claim~\ref{cl: step2} $(iii)$,
and $|R_e|\geq \frac12(\log\Delta)^8$ by Claim~\ref{cl: reserve} $(i)$.
We thus obtain a proper edge colouring $c$ of $G$ from the given lists.
As in this final step for every edge $uv\in E$ with $d(u)=d(v)$ we have only assigned colours to its adjacent edges in $U'_u\cup U'_v$, the palettes of $u$ and $v$ remained distinct since by Claim~\ref{c2: step3} we had
$|S'_u \triangle S'_v| \geq \frac1{32} (\log\Delta)^2$ and
$|U'_u\cup U'_v| \leq \frac12\cdot\frac1{32} (\log\Delta)^2$,
and consequently now $|S_c(u) \triangle S_c(v)| \geq \frac{1}{2}\cdot\frac1{32} (\log\Delta)^2$.
This finishes our construction, as the remaining adjacent vertices (with distinct degrees) are set distinguished by definition.
$\Box$

\section{Concluding remarks}

Note that we have in fact showed above the following stronger thesis than in Theorem~\ref{main_theorem_list_Hatami_auxiliary}.
\begin{theorem}\label{main_theorem_list_Hatami_auxiliary_enhanced}
There is a constant $C_0$ and $\Delta_0$ such that for every multigraph $G=(V,E)$ with edge multiplicity at most $2$, maximum degree $\Delta\geq\Delta_0$ and minimum degree $\delta\geq\frac{\Delta}{4}$, and any set $\{L_e\}_{e\in E}$ of lists of cardinalities $\Delta + \lfloor C_0\Delta^\frac{1}{2}(\log\Delta)^4\rfloor$,
there exists a choice of colours from these lists yielding a proper edge colouring $c$ of $G$ such that $|S_c(u) \triangle S_c(v)| \geq \frac{1}{2}\cdot\frac1{32} (\log\Delta)^2$ for every edge $uv\in E$ with $d(u)=d(v)$.
\end{theorem}
This almost immediately implies the following total correspondent of Theorem~\ref{main_theorem_list_Hatami} for simple graphs,
where the graph invariant ${\rm ch}''_a(G)$ is defined the same as ${\rm ch}'_a(G)$ but for total, not edge
colourings $c$ and with the palette of a vertex $v$ defined as $S'_c(v):=c(v)\cup\{c(e):e\in E(v)\}$.
\begin{corollary}\label{main_theorem_list_Hatami_total}
There is a constant $C'$ such that
$${\rm ch}''_a(G)\leq \Delta + C'\Delta^\frac{1}{2}(\log\Delta)^4$$
for every graph $G=(V,E)$ with maximum degree $\Delta$.
\end{corollary}

\begin{pf}
Analogously as previously it is sufficient to show the thesis for graphs with maximum degree $\Delta$ large enough, e.g. for $\Delta\geq\max\{\Delta_0,e^{16}\}$ and with $C'=C_0+1$ where $C_0$ and $\Delta_0$ are the constants from Theorem~\ref{main_theorem_list_Hatami_auxiliary_enhanced}.
Fix any set of lists $\{L_x\}_{x\in V\cup E}$ of lengths $\Delta + \lfloor C'\Delta^\frac{1}{2}(\log\Delta)^4\rfloor$.
First we greedily choose an auxiliary proper vertex colouring $c$ of $G$ from the given lists and for every edge $uv\in E$ we remove $c(u)$ and $c(v)$ from the list of $uv$. Denote the obtained set of edge lists by $\{L'_e\}_{e\in E}$ and note that $|L'_e|\geq \Delta + \lfloor C_0\Delta^\frac{1}{2}(\log\Delta)^4\rfloor$.
Then as long as the minimum degree is less than $\Delta/4$ we iteratively take two copies of a currently analyzed graph and add edges between the corresponding vertices of degree less than $\Delta/4$, until finally we obtain a graph $G'$ with minimum degree $\delta\geq \Delta/4$.
Fix any list $L'_e$ with $|L'_e|= \Delta + \lfloor C_0\Delta^\frac{1}{2}(\log\Delta)^4\rfloor$ for every edge $e$ of $G'$ without a list assigned.
By Theorem~\ref{main_theorem_list_Hatami_auxiliary_enhanced} there is a proper edge colouring $c'$ of $G'$ from the lists $L'_e$ such that $|S_{c'}(u) \triangle S_{c'}(v)| \geq \frac{1}{2}\cdot\frac1{32} (\log\Delta)^2\geq 4$ for every edge $uv$ with equal degrees of $u$ and $v$.
This restricted to the edges of the original graph $G$ retains this property for the adjacent vertices $u,v$ in $G$ with $d(u)=d(v)\geq\Delta/4$.
We then colour every vertex $v$ with $c(v)$ and next greedily change the colour of every vertex $v$ with $d(v)<\Delta/4$
so that the obtained total colouring $c''$ of $G$ is proper and there is no set conflict between any vertex $v$ with $d(v)<\Delta/4$
and its neighbours in $G$. At the same time, by our construction,
$|S'_{c''}(u) \triangle S'_{c''}(v)| \geq 2$ for every edge $uv$ with $d(u)=d(v)\geq\Delta/4$ in $G$.
$\Box$
\end{pf}

\end{document}